\documentclass{amsart}
\newcommand{\note}{\noindent {\bf Notation. }}

\newcommand{\ws}{\hspace{4pt}}
\newtheorem{theorem}{Theorem}

\newtheorem{remark}{Remark}
\newtheorem{proposition}{Proposition}

\newtheorem{lemma}{Lemma}

\usepackage[]{amssymb,amsopn}
\usepackage[]{amsmath}
\usepackage[]{eucal}
\usepackage[]{latexsym}
\usepackage{color}

\makeatletter
\@namedef{subjclassname@2020}{%
  \textup{2020} Mathematics Subject Classification}
\makeatother

\begin{document}

\title[Bessel trace inequality]{Trace inequality with Bessel convolution}
\author{Mouna Chegaar and \'A. P. Horv\'ath}

\thanks{The second author is supported by the NKFIH-OTKA Grant K132097.}

\begin{abstract} Considering potentials defined by Bessel kernel with Bessel convolution a Kerman-Sawyer type characterization of trace inequality is given. As an application an estimate on the least eigenvalue of Schr\"odinger-Bessel operators is derived.

\end{abstract}
\maketitle

  \section{introduction}

 After  several authors we call trace inequality the next type of inequalities
 $$\|Jf\|_{p,\mu} \le C \|f\|_{p,\nu},$$
 where $J$ is a certain integral operator, $\mu$ and $\nu$ are positive measures. The reason is that $\mu$ can be the restriction of $\nu$ to a subspace see \cite {ks} and \cite{st}, or considering the Sobolev type version it generalizes $\|u\|_{L^q(\Omega,d\mu)}\le C\|\nabla_m u\|_{L^p(\Omega,d\nu)}$, see \cite{v}.
  The notion is strongly related to non-linear potential theory and has several different applications. For instance it can be applied to characterize positive measures which are multipliers for a pair of potential spaces, embedding theorems for Sobolev spaces, see e.g. \cite{mv}, \cite{ms}, \cite{lg}, to obtain bounds
for the Hardy-Littlewood maximal function, see \cite{op}. It is also useful in applications to criteria of solvability for non-linear partial differential equations, see e.g. \cite{hmv}, \cite{np}, \cite{ap}, and for estimation of eigenvalues of partial differential operators, see e.g. \cite{ks}. \\
 Some connections with nonlinear capacity in case of general radially decreasing convolution kernels are described in \cite{ah}. In the celebrated paper of Maz'ya and Verbitsky (see \cite{mv}) the equivalence of the trace inequality to capacitary and potential-type inequalities is proved in  Riesz and Bessel cases. In Bessel case, some additions  and additional results can be found in \cite{o}.\\
 As potential theory developed by Gauss in order to describe gravitational and electrostatic fields, it is ball-symmetric and is strongly connected to the Laplace operator. So partial differential operators arise in the applications are connected to the Laplace operator as well. Below we change the symmetry of the space applying Bessel translation and convolution rather than the standard one. This implies that instead of Laplace, the Laplace-Bessel operator (see \eqref{lbo}) comes into the picture which is also an important tool in physics. Moreover it results investigations in weighted spaces and so the trace inequality takes on the general form above. (For weighted potential see e.g. \cite{a}, for potential with Bessel convolution see \cite{h}.) In Theorem \ref{vm} we show the equivalence of the trace inequality with a "cube condition", c.f. \cite{ks}. Since it is much simpler to check a condition on cubes than on general compact sets, it proved to be a useful tool for applications, see \cite{ks} and Theorem \ref{eiv} below.

 The paper is organized as follows. The following section contains the main results and the definitions that are essential for their formulation. In the third section we collected all the properties of the defined notions which are used later. Moreover this section contains the proofs of some technical lemmas. The fourth section is devoted to the proof of Theorem \ref{vm}. As an application of Theorem \ref{vm} in the last section we deal with the Schr\"odinger-Bessel operator.

  \section{Notation, main results}
  Let $ \mathbb{R}^{n}_{+}:= \{x = (x_{1},\dots ,x_{n}): x_{i}\geq 0, \  i = 1, \dots ,n.\}$. $B(x,r)$ and $Q(x,r)$ stand for the positive part of balls and cubes centered at $x\in \mathbb{R}^n_+$ and of radius (half-side) $r$ (i.e. e.g. $Q(x,r)=\times_{i=1}^n (x_i-r,x_i+r)\cap \mathbb{R}^{n}_{+}$). $\lambda$ is the n-dimensional Lebesgue measure.
  $a = a_{1}, \dots , a_{n}$ is a multiindex. Let $\mathrm{ E} \subset\mathbb{R}^{n}_{+} $ and $ \mathcal{M}(E)$ stands for Radon
measures supported on $\mathrm{E}$. If $\mu$ is a measure and $K$ is a set, then $\mu_{|_K}$ is the restriction of $\mu$ to $K$. If $\mu\in\mathcal{M}(E)$ for some $E,$ $ d\mu_{a}(x):=x^{a} d\mu (x).$ A measure $\mu$ is doubling if there is a constant $C$ such that for all $r>0$ and $x\in \mathbb{R}^{n}_{+}$ $\mu(B(x,2r))\le C\mu(B(x,r))$.

Define the Banach space  $ L^{p}_{a}$ as follows.
\begin{equation}\label{}
  \|f \|^{p}_{p,a}=\int_{\mathbb{R}^{n}_+} |f(x)|^{p} d\lambda_{a}(x).
\end{equation}
  and as usual
  \begin{equation}
    L^{p}_{a}:= L^{p}_{a}(\mathbb{R}^{n}_{+})= \{f : \|f \|_{p,a} < \infty \},  \ \ \    \mathrm{L}^{p+}_{a}:= \{f \in L^{p} : f \geq 0 \}.
  \end{equation}
  The dual index $p$ is defined by  $ \frac{1}{p}+\frac{1}{p'}= 1.$

 \subsection{Bessel Translation}
   Let $a_i= 2\alpha_{i} + 1, \ \alpha_{i} > -\frac{1}{2}$, $i = 1, \dots , n$,\\ $|a|= \sum_{i=1}^{n}(2 \alpha_{i} + 1)$, $t \in \mathbb{R}^{n}_{+}$.

The Bessel translation of a function, $f$ (see e.g. [14], [16], [20]) is
\begin{equation}
  T^{t}_{a}f(x)= T^{t_{n}}_{a_{n}} \dots T^{t_{1}}_{a_{1}} f (x_{1}, \dots , x_{n}),
\end{equation}
where

$$ T^{t_{i}}_{a_{i}} f (x_{1}, \dots , x_{n})=\frac{\Gamma(\alpha_{i}+1)}{\sqrt{\pi} \Gamma(\alpha_{i}+\frac{1}{2})}$$ \begin{equation}\label{a} \times\int_{0}^{\pi} f(x_{1},\dots,\sqrt{x_{i}^{2}+ t_{i}^{2} - 2 x_{i}^{2} t_{i}^{2}\cos\vartheta_{i}} , x_{i+1},\dots,x_{n})\sin^{2 \alpha_{i}} \vartheta_{i}d\vartheta_{i}.
\end{equation}

The generalized convolution with respect to Bessel translation is
\begin{equation}
  f\ast_{a}g=\int_{\mathbb{R}^{n}_{+}}  T^{t}f(x)g(x) d \lambda _{a}(x).
\end{equation}

\subsection{Radially decreasing and Bessel kernels}

Let $g$ be a non-negative lower semi-continuous, non-increasing function on $\mathbb{R}_{+} $
for which
\begin{equation}\label{f}
  \int_{0}^{1} g(t)t^{n+|a|-1} dt < \infty .
\end{equation}
Then  $\kappa:=g(|x|)$ is a radially decreasing kernel on $ \mathbb{R}^{n} $.\\

\medskip

The modified Bessel function of the second kind, $K_{\alpha}$ is defined as follows.
\begin{equation}\label{g}
  i^{-\alpha}J_{\alpha}(ix)=\sum_{k=0}^{\infty}\frac{1}{k!\Gamma(k+\alpha+1)}\left(\frac{x}{2}\right)^{2k+\alpha},
\end{equation}
where $J_{\alpha}$ is the Bessel function, and
\begin{equation}\label{h}
  K_{\alpha}(x)=\frac{\pi}{2}\frac{i^{\alpha}J_{-\alpha}(ix)-i^{-\alpha}J_{\alpha}(ix)}{\sin \alpha \pi}.
\end{equation}

The Bessel kernel is
\begin{equation}\label{k}
  G_{a,\nu}(x):=\frac{2^{\frac{n-a-\nu}{2}+1}}{\Gamma(\frac{\nu}{2})\prod_{i=1}^{n}\Gamma (\alpha_{i}+1)}\frac{K_{\frac{n+|a|-\nu}{2}}(|x|)}{|x|^{\frac{n+|a|-\nu}{2}}}.
\end{equation}

\subsection{Schr\"odinger-Bessel differential operator}

Schr\"odinger-Bessel operator is
\begin{equation}\label{sre}H_a:=-\Delta_a-v,\end{equation}
where $\Delta_a$ is the Laplace-Bessel operator, i.e.
\begin{equation}\label{lbo}\Delta_a u=\sum_{i=1}^n u_{x_ix_i}+\frac{2\alpha_i+1}{x_i}u_{x_i},\end{equation}
and $v$ is nonnegative and locally integrable.

\subsection{Results}

The next result gives a Kerman-Sawyer type characterization of trace inequality. It is proved for Riesz and Bessel kernels with respect to the standard convolution, see \cite{ks}, \cite{mv}, \cite{o}. Here we state a similar result with Bessel convolution.

\medskip

\begin{theorem}\label{vm} Let $1<p<\infty$, $0<\nu$. Let $\gamma \in  \mathcal{M}(\mathbb{R}^{n}_+))$ be a doubling measure. Then the following properties of $\gamma$ are equivalent.

\noindent There is a constant $A_{1}$ such that for all $ f \in L_{a}^p$
\begin{equation}\label{I}
  \left(\int_{\mathbb{R}^{n}_{+}} | G_{a,\nu} *_{a}f|^{p} d\gamma_a \right)^{\frac{1}{p}} \leq A_{1} \|f\|_{p,a}.
\end{equation}

 \noindent   There is a constant $A_{2}$ such that for all cubes $Q$
\begin{equation}\label{Q} \int_Q ((G_{a,\nu}*_a\gamma|_Q)(x))^{p'}x^adx \le A_2\gamma_a(Q),\end{equation}

 \noindent  The potential $G_{a,\nu}*_a\gamma$ is finite a. e. and there is a constant $A_{3}$ such that
\begin{equation}\label{V}G_{a,\nu}*_a(G_{a,\nu}*_a\gamma)^{p'}(x)\le A_3 \ (G_{a,\nu}*_a\gamma)(x). \end{equation}
Moreover the least possible constants $A_i$ are all proportional with each other with constants that are independent of $\gamma$.
\end{theorem}

\medskip

The class of measures characterized by \eqref{V} turned out to be useful in applications to nonlinear equations, see e.g. \cite{hmv}. Below we apply characterization \eqref{Q} to get some eigenvalue estimate to a Schr\"odinger-type operator.

Let $l(Q)$ be the edge lengths of the cube $Q$, and $\mu(Q)$ is the volume of $Q$ with respect to the measure $\mu$.

\medskip

\begin{theorem}\label{eiv} Let $H_a$ be as defined by \eqref{sre}. Assume further that $v(x)$ defines a doubling measure. Then there are two positive constants $A$ and $B$ depending only on $n$ and $a$ such that the least eigenvalue $\lambda_1$ of $H_a$ satisfies
$$L\le -\lambda_1 \le U,$$
where
$$L= \sup\left\{l(Q)^{-2}:\frac{1}{v_a(Q)}\int_Q G_{a,2}*_av|_Q(x)x^av(x)dx\ge B\right\};$$ $$ U= \sup\left\{l(Q)^{-2}:\frac{1}{v_a(Q)}\int_Q G_{a,2}*_av|_Q(x)x^av(x)dx\ge A\right\}.$$\end{theorem}

\section{Preliminaries}

In this section some useful properties of the above defined notions are collected and derived.

\subsection{Bessel translation}

Bessel translation can also be expressed as an integral with respect to a kernel function:
\begin{equation}\label{b}
  T^{t_{i}}_{a_{i}} f (x_{1}, \dots , x_{n})=\int_{0}^{\infty} K(x_{i},t_{i},z_{i})f(z_{1},\dots,z_{n})d\lambda_{a_{i}}(z_{i}),
\end{equation}
where for any $ \alpha > -\frac{1}{2} ,x,t\geq 0 $
\begin{equation}\label{c}
  K(x,t,z)=\left\{
             \begin{array}{ll}
               \frac{\pi^{\alpha+\frac{1}{2}} \Gamma(\alpha+1)}{2^{2\alpha-1}\Gamma(\alpha+\frac{1}{2})} \frac{[((x+t)^{2}-z^{2})(z^{2}-(x-t)^{2})]^{\alpha-\frac{1}{2}}}{(xtz)^{2\alpha}}, \ws  |x-t|< z <|x+t|; \\
               0, \ws \mbox{otherwise.}
             \end{array}
           \right.
\end{equation}

Since $a$ is fixed subsequently we omit it from the notation of the translation.\\
The formula above shows that Bessel translation is positive operator. It is also symmetric and bounded, that is

\begin{equation}
  T^{t}f(x) = T^{x} f(t).
\end{equation}

\begin{equation}\label{d}
  \| T^{t}f(x) \|_{p,a} \leq \|f\|_{p,a}, \ \ \ 1 \le p \le  \infty ,
\end{equation}
see e.g. \cite{le}.

Subsequently we need the estimations below.

\begin{lemma}\label{tk}\cite{h} $\mathrm{supp}T^t\chi_{B_+(0,r)}(x)=\overline{B_+(x,r)}$, $\mathrm{supp}T^t\chi_{[0,r)^n}(x)=\times_{i=1}^n[x_i-r,x_i+r]_+=:T_+(x,r)$.
There is a $c>0$ such that for all $x\in \mathbb{R}^n_+$, $t\in B_+(x,r)$
\begin{equation}\label{kf}T^t\chi_{B_+(0,r)}(x)\le c \prod_{i=1}^n\min\left\{1,\left(\frac{r}{x_i}\right)^{a_i}\right\}.\end{equation}
There is a $c>0$ such that for all $x\in \mathbb{R}^n_+$, $t\in T_+\left(x,\frac{r}{2}\right)$
\begin{equation}\label{kf1}T^t\chi_{[0,r)^n}(x)\ge c \prod_{i=1}^n\min\left\{1,\left(\frac{r}{x_i}\right)^{a_i}\right\}.\end{equation}
\end{lemma}

To formulate the next technical lemma we introduce the abbreviations:

$$d(x,t,\vartheta):=\sum_{i=1}^n x_i^2+t_i^2-2x_it_i\cos\vartheta_i,$$
$$D:=\left|\sqrt{d(t,u,\vartheta)}-\sqrt{d(x,y,\vartheta)}\right|, \ws \ws N:=d(t,u,\vartheta)+d(x,y,\vartheta).$$

\medskip

\begin{lemma} Let $\vartheta \in [0,\pi]^n$ be fixed. Then we have
\begin{equation}\label{v5}D \le 2 \max\{\|x-t\|, \|y-u\|\}\left(3+\sqrt{2}\frac{\|y-t\|}{\sqrt{N}}\right).\end{equation}
\begin{equation}\label{v51}d(z,y,\vartheta)\le C d(x,y,\vartheta) \ws \ws \mbox{if} \ws \|z-y\|\le B \|x-y\|.\end{equation}
\end{lemma}

\proof First we show \eqref{v5}.
$$D \le \frac{\sum_{i=1}^n(t_i-x_i)(t_i+x_i-2y_i\cos\vartheta_i)+(u_i-y_i)(u_i+y_i-2t_i\cos\vartheta_i)}{\sqrt{d(t,u,\vartheta)+d(x,y,\vartheta)}}$$
$$\le \sqrt{2}\max\{\|x-t\|, \|y-u\|\}\sqrt{\frac{N +\sum_{i=1}^nH_i}{N}},$$
where $H_i=4\cos^2\vartheta_i(y_i^2+t_i^2)+2(t_ix_i+u_iy_i)-8t_iy_i\cos\vartheta_i-2t_iu_i\cos\vartheta_i-2x_iy_i\cos\vartheta_i$. Thus
$$D \le  \sqrt{2}\max\{\|x-t\|, \|y-u\|\}\sqrt{\frac{2N+\sum_{i=1}^nS_i}{N}},$$
where
$$S_i:=4(\cos^2\vartheta_i(y_i^2+t_i^2)-2\cos\vartheta_it_iy_i).$$  Let $N_i$ be the $i$th term in the denominator. If $S_i\le 0$, then we estimate it by 0. If $\cos\theta_i< 0$ ($S_i>0$), then $S_i\le 2\times 4(y_i^2+t_i^2)\le 8 N_i$. If $S_i>0$ and $\cos\theta_i> 0$
 then $\cos\theta_i> \frac{2t_iy_i}{t_i^2+y_i^2}$ and so $S_i$ is increasing in $\cos\theta_i$ i.e. $S_i\le 4(t_i-y_i)^2$.  Thus
$$D\le  \sqrt{2}\max\{\|x-t\|, \|y-u\|\}\left(\sqrt{\frac{2N}{N}}+\sqrt{\frac{\sum_{ \cos\theta_i< 0} S_i}{N}}+\sqrt{\frac{\sum_{ \cos\theta_i> 0} S_i}{N}}\right)$$ $$\le  \sqrt{2}\max\{\|x-t\|, \|y-u\|\}\left(\sqrt{2}+\sqrt{8}+2\frac{\|y-t\|}{\sqrt{N}}\right).$$

To prove \eqref{v51} observe that
$$d(z,y,\vartheta)= \|z-y\|^2+4\sum_{i=1}^ny_iz_i\sin^2\frac{\vartheta_i}{2}\le B^2\|x-y\|^2+4\sum_{i=1}^ny_iz_i\sin^2\frac{\vartheta_i}{2}.$$
Let $A>1$. If $Ax_i>z_i$ for some index $i$, then
\begin{equation}\label{A1} 4y_iz_i\sin^2\frac{\vartheta_i}{2}\le A 4y_ix_i\sin^2\frac{\vartheta_i}{2}.\end{equation}
If $Ax_i\le z_i$ for some $i$, then $(A-1)x_i\le z_i-x_i$. In this case we take
$$y_iz_i=x_iy_i+(z_i-x_i)(y_i-x_i)+(z_i-y_i)x_i+(y_i-x_i)^2+(y_i-x_i)x_i$$ $$\le x_iy_i+(z_i-x_i)(y_i-x_i)+(z_i-y_i)\frac{1}{A-1}((z_i-y_i)+(y_i-x_i))+(y_i-x_i)^2$$ $$+(y_i-x_i)\frac{1}{A-1}((z_i-y_i)+(y_i-x_i)).$$
Thus
$$d(z,y,\vartheta)\le B^2\|x-y\|^2+4\sum_{1\le i\le n \atop Ax_i\le z_i}y_iz_i\sin^2\frac{\vartheta_i}{2}+4\sum_{1\le i\le n \atop Ax_i> z_i}(\cdot)$$ $$\le B^2\|x-y\|^2+4\sum_{1\le i\le n \atop Ax_i\le z_i}y_ix_i\sin^2\frac{\vartheta_i}{2}+4\|z-x\|\|y-x\|+\frac{1}{A-1}\|z-y\|^2$$ $$+\frac{1}{A-1}\|z-y\|\|x-y\|+\|x-y\|^2+\frac{1}{A-1}\|z-y\|\|x-y\|+\frac{1}{A-1}\|x-y\|^2 $$ $$+4 A\sum_{1\le i\le n \atop Ax_i> z_i}y_ix_i\sin^2\frac{\vartheta_i}{2}\le  (A+1)4\sum_{i=1}^ny_ix_i\sin^2\frac{\vartheta_i}{2} + \|x-y\|^2\left(2B+\frac{(B+1)^2}{A-1}\right),$$
i.e. \eqref{v51} is proved by $C=\max\left\{A+1, 2B+\frac{(B+1)^2}{A-1}\right\}$.

\medskip

\begin{remark}

In particular
$$ D\le 10M, \ws \mbox{if} \ws x=t , \ws  \|u-y\|<M.$$
Indeed, $\sqrt{N} \ge\frac{1}{\sqrt{2}}(\|y-x\|+\|u-t\|)$ and $\|y-t\|=\|y-x\|$.
\end{remark}

\medskip

\subsection{Bessel convolution}

The convolution defined by Bessel translation similarly to the standard one, is symmetric and it fulfils the Young's inequality, i.e.

\begin{equation}
   f\ast_{a}g =  g\ast_{a}f,
\end{equation}
 and if  $1\leq p, q, r\leq \infty $  with  $ \frac{1}{r}=\frac{1}{p}+\frac{1}{q}-1 $ ;   if \  $f\in L^{p}_{a}$ and
$ g \in  L^{q}_{a},$ then

\begin{equation}\label{e}
  \| f\ast_{a}g \|_{r,a} \leq \|f\|_{p,a}\|g\|_{q,a},
\end{equation}
see \cite[(3.178)]{ss}.

\medskip

\subsection{Radially decreasing kernels with Bessel convolution}
The next lemma with general radially decreasing kernel is useful tool of the investigation below.

\medskip

\begin{lemma}\label{l1} Let $\kappa(t)=g(|t|)$ be a radially decreasing kernel. Then
$$\kappa *_a\mu(x)= \int_0^\infty \chi_{B(0,r)}*_a\mu(x)(- g'(r))dr.$$
\end{lemma}

\proof Let $R>0$ be arbitrary. Define the next measure: $d\Theta_{x,a}(t)=t^adT^x\mu(t)$. Since $g$ is differentiable a.e. we have
$$\int_0^R \chi_{B(0,r)}*_a\mu(x)(- g'(r))dr=\int_0^R(- g'(r))\int_{B(0,r)}1d\Theta_{x,a}(t)dr$$ $$=\int_{B(0,r)}\int_{|t|}^R(- g'(r))drd\Theta_{x,a}(t)$$ $$=\int_{B(0,r)}g(|t|)d\Theta_{x,a}(t)-g(R)\int_{B(0,r)}1d\Theta_{x,a}(t).$$
Then commutativity of convolution implies that
$$\int_0^R \chi_{B(0,r)}*_a\mu(x)(- g'(r))dr=(\kappa\chi_{B(0,R)})*_a\mu(x) - g(R)(\chi_{B(0,R)}*_a\mu(x)) .$$
Similarly
$$(\kappa \chi_{\mathbb{R}^n_+\setminus B(0,R)})*_a\mu(x)=\int_R^\infty \chi_{B(0,r)}*_a\mu(x)(- g'(r))dr- g(R)(\chi_{B(0,R)}*_a\mu(x)) .$$
Added up the result follows.

\medskip

In the main part of the paper we restrict the investigation to Bessel kernel.  It has the next asymptotic properties, see e.g. \cite[page 12]{ah}. \\
Around zero we have

 \begin{equation}\label{i}
     K_{\alpha}(r)\sim  \left\{
                         \begin{array}{ll}
                            -ln \frac{r}{2}-c, & \text {if}  \   \alpha=0;   \\
                          C(\alpha)r^{-\alpha}, & \text{if}   \  \alpha > 0,
                         \end{array}
                       \right.  \end{equation}
and around infinity
 \begin{equation}\label{j}
    K_{\alpha}(r)\sim \frac{c}{\sqrt{r}}e^{-r}.
 \end{equation}

 Recall the definition of Bessel kernel, \eqref{k}.  $G_{a,\nu}$ has the next semi-group property (see \cite[Lemma 4.3]{sek})
 \begin{equation}\label{sgp} G_{a,\nu}*_a G_{a,\mu}= G_{a,\nu+\mu}, \ws \ws \nu, \mu >0.\end{equation}

 \subsection{Hankel transformation}
 Let
 $$j_\alpha(x)=\frac{2^\alpha \Gamma(\alpha+1)}{x^\alpha}J_\alpha(x)$$
 be the entire Bessel function. Then the Hankel transform of a function $f\in L^1_a(\mathbb{R}_{+}^n)$ is
 $$\mathcal{H}_af(\xi)=\hat{f}(\xi)=\int_{\mathbb{R}_{+}^n}f(x)\prod_{i=1}^nj_{\alpha_i}(x_i\xi_i)x^adx.$$
 Below we need the next properties of Hankel transformation.
 \begin{equation}\label{hkonv} \mathcal{H}_af*_ag=\hat{f}\hat{g},\end{equation}
 \begin{equation}\label{hb}\mathcal{H}_a(-\Delta_a u)(\xi)=|\xi|^2\hat{u}(\xi)\end{equation}
 and denoting by
 $$\langle f,g \rangle_a:=\int_{\mathbb{R}_{+}^n}f(x)g(x)x^adx,$$
 \begin{equation}\label{hn}\langle f,g \rangle_a= \frac{2^{n-|a|}}{\prod_{i=1}^n\Gamma^2(\alpha_i+1)}\langle\hat{f},\hat{g}\rangle_a ,\end{equation}
 see e.g. \cite[(1.83), (1.95)]{ss}.
 Finally for $\nu>0$ we have (see \cite[Lemma 4.1]{sek}) that
 \begin{equation}\label{hG}  G_{a,\nu}(x)=\mathcal{H}_a^{-1}\frac{1}{(1+|\xi|^2)^{\frac{\nu}{2}}}.\end{equation}

\section{Proof of Theorem \ref{vm}.}

\subsection{\eqref{I}$\Rightarrow$ \eqref{Q}}

The first implication in Theorem \ref{vm} can be stated in more general form. For the corresponding statement with respect to standard convolution we refer to \cite[Theorem 7.2.1]{ah}.

 \medskip

 \begin{proposition}
   Let $\kappa \in L^1_a$ be a radially decreasing convolution kernel, and let $\mu \in  \mathcal{M}(\mathbb{R}^{n}_+)$ and  $ 1 < p \le q < \infty $. Then if
there is a constant $B_{1}$ such that for all $ f \in L_{a}^p$
\begin{equation}\label{l}
  \left(\int_{\mathbb{R}^{n}_{+}} | \kappa \ast_{a}f|^{q} d\mu_a \right)^{\frac{1}{q}} \leq B_{1} \|f\|_{p,a},
\end{equation}
Then there is a constant $B_{2}$ such that for all compact sets $K$
\begin{equation}\label{m}
  \|\kappa \ast_{a}\mu|_{K} \|_{p',a}\leq B_{2} \mu_a (K)^{\frac{1}{q'}}.
\end{equation}
 \end{proposition}

\proof

\eqref{l} $\Rightarrow$ \eqref{m}. Let $f \in L_{a}^p$. This implication is immediate from
  $$\int_{\mathbb{R}^{n}_{+}} f(x) (\kappa \ast_{a}\mu|_{K})(x)x^a dx  = \int_{\mathbb{R}^{n}_{+}} f(x) \int_{\mathbb{R}^{n}_{+}}  T^{t}\kappa (x) t^{a} d\mu_{|_K}(t) x^{a} dx $$
 $$= \int_{\mathbb{R}^{n}_+} t^{a} \int_{\mathbb{R}^{n}_{+}} f(x) T^{t} \kappa (x) x^{a}dx d\mu_{|_K}(t) = \int_{\mathbb{R}^{n}_{+}} (\kappa \ast_{a} f )(t) t^{a} d\mu_{|_K}(t) $$ $$= \int_{K} (\kappa \ast_{a} f ) d\mu_a(t).$$
 Applying H\" older's inequality then combining with \eqref{l} and \eqref{d} yields
$$\left|\int_{\mathbb{R}^{n}_{+}} f(x) (\kappa \ast_{a}\mu|_{K})(x)x^a dx \right| \leq \| \kappa \ast_{a} f  \|_{q, \mu_a}\mu_a(K)^{\frac{1}{q'}} \leq B_{1} \|f\|_{p,a} \ \mu_a(K)^{\frac{1}{q'}},$$
 for all   $f \in L^p_a.$
 Finally we have
 $$\|\kappa \ast_{a}\mu|_{K} \|_{p',a}  \leq  B_{2}  \mu_a (K)^{\frac{1}{q'}}.$$

\subsection{\eqref{Q}$\Rightarrow$ \eqref{V}}
To prove that \eqref{Q} implies \eqref{V}, first let us observe that \eqref{m} implies \eqref{Q}, and assume now that $\gamma$ is doubling. We need some further observations.

\medskip

\begin{lemma} If $\gamma$ is doubling, then $\gamma_a$ is doubling as well and
\begin{equation}\label{v9}(\chi_{B(0,r)}*_a\gamma)(x) \sim r^{|a|}\gamma(B(x,r)).\end{equation}\end{lemma}

\proof

The first statement is obvious. To see the second one notice that if $t\in B(x,r)$, then $\min\{1,\frac{r}{x_i}\}t_i\le 2r$. By Lemma \ref{tk} we have
$$(\chi_{B(0,r)}*_a\gamma)(x)\le c\int_{B(x,r)}r^{|a|}d\gamma(t),$$
which gives the upper estimate (independently of doubling condition). On the other hand
$$(\chi_{B(0,r)}*_a\gamma)(x)\ge (\chi_{Q(0,\frac{r}{\sqrt{n}})}*_a\gamma)(x)\ge c\int_{Q(x,\frac{r}{2\sqrt{n}})}\prod_{i=1}^n\min\left\{1,\left(\frac{r}{x_i}\right)^{a_i}\right\}t^ad\gamma(t)$$ $$\ge c\int_{Q(x+\frac{r}{4}e,\frac{r}{4\sqrt{n}})}\prod_{i=1}^n\min\left\{1,\left(\frac{r}{x_i}\right)^{a_i}\right\}t^ad\gamma(t)\ge c r^{|a|}B\left(x+\frac{r}{4}e,\frac{r}{4\sqrt{n}}\right)$$ $$\ge cr^{|a|}\gamma\left(B\left(x+\frac{r}{4}e,2r\sqrt{n}\right)\right)\ge cr^{|a|}\gamma(B(x,r)),$$
where $e=(1,1,\dots ,1)$ and the doubling property is used.

\medskip

\begin{remark}\label{rem2} If $\gamma$ is doubling, then
$$\gamma(Q(x,r))\prod_{i=1}^n\left(\max\{x_i,r\}\right)^{a_i} \sim \gamma_a(Q(x,r)).$$
\end{remark}

Indeed,
$$\gamma_a(Q(x,r))=\int_{Q(x,r)}t^ad\gamma(t)\le \prod_{i=1}^n(x_i+r\sqrt{n})^{a_i}\gamma(Q(x,r)),$$
and on the other hand
$$\gamma(Q(x,r))\le \gamma(Q(x+r\sqrt{n}e,2r))\le c\gamma\left(Q\left(x+r\sqrt{n}e,\frac{r}{2}\right)\right)$$ $$\le c\frac{1}{\left(x+\frac{r}{2}\sqrt{n}e\right)^a}\int_{Q\left(x+r\sqrt{n}e,\frac{r}{2}\right)}t^ad\gamma(t)\le c\frac{1}{\left(x+\frac{r}{2}\sqrt{n}e\right)^a}\gamma_a\left(Q\left(x,\frac{3}{2}r\right)\right)$$ $$\le c\frac{1}{\left(x+\frac{r}{2}\sqrt{n}e\right)^a}\gamma_a(Q(x,r)).$$

\medskip

\begin{lemma}\label{v4} Let $\gamma \in \mathcal{M}$ such that it satisfies \eqref{Q}. Let $c_0$ and $A>1$ be arbitrary but fixed constants. Then
\begin{equation}\label{v41}\gamma_a(Q(x,r))\le c  r^{n-p\nu}\prod_{i=1}^n\left(\max\{x_i,r\}\right)^{a_i},\end{equation}
if $r\le c_0$.\\
If $r>c_0$, then
\begin{equation}\label{v42}\gamma_a(Q(x,r))\le cr^{-n(p-1)-|a|p}\left(g\left(\frac{r}{A}\right)\right)^{-p}\left(\max\{x_i,r\}\right)^{a_i},\end{equation}
where $c$ depends on $c_0$ and $A$.
\end{lemma}

\proof  To prove \eqref{v41} let $z\in Q(x,r)$. Considering that if $\varrho>2\sqrt{n}r$\\ then $Q(x,r)\subset B(z,\varrho)$ and we have
$$(G_{a,\nu}*_a\gamma|_{Q(x,r)})(z)= \int_0^\infty -g'(\varrho)\int_{B(z,\varrho)\cap Q(x,r)}T^u\chi_{B(0,\varrho)}(z)u^ad\gamma(u)d\varrho$$
$$\ge c\int_{2\sqrt{n}r}^\infty -g'(\varrho)\int_{Q(x,r)}T^u\chi_{B(0,\varrho)}(z)u^a d\gamma(u)d\varrho $$ $$\ge c\int_{2\sqrt{n}r}^\infty-g'(\varrho)\int_{Q(x,r)}T^u\chi_{Q\left(0,\frac{\varrho}{\sqrt{n}}\right)}(z)u^a d\gamma(u) d\varrho.$$
If $\varrho>4\sqrt{n}r$ then $Q(x,r)\subset Q\left(0,\frac{\varrho}{2\sqrt{n}}\right)$, and by Lemma \ref{tk}\\ $T^u\chi_{Q\left(0,\frac{\varrho}{\sqrt{n}}\right)}(z)\ge c \prod_{i=1}^n\min\left\{1,\left(\frac{\varrho}{z_i}\right)^{a_i}\right\}$.\\ Let us suppose first that $x_i>r$ $i=1, \dots , n$. Then $\frac{\varrho}{z_i}\ge \frac{\varrho}{x_i+r}\ge \frac{\varrho}{2x_i}$. Thus
$$(G_{a,\nu}*_a\gamma|_{Q(x,r)})(z)\ge c\int_{4\sqrt{n}r}^\infty -g'(\varrho)\int_{Q(x,r)}\prod_{i=1}^n\min\left\{1,\left(\frac{\varrho}{x_i}\right)^{a_i}\right\}u^a d\gamma(u)d\varrho$$ $$ \ge c g(4\sqrt{n}r)\prod_{i=1}^n\min\left\{1,\left(\frac{r}{x_i}\right)^{a_i}\right\}\gamma_a(Q(x,r)).$$
Now integrating over $Q$ with respect to $x^a$ by the assumption (c.f. \eqref{Q}) we have
$$c x^ar^n g(4\sqrt{n}r)^{p'}\left(\prod_{i=1}^n\min\left\{1,\left(\frac{r}{x_i}\right)^{a_i}\right\}\right)^{p'}\gamma_a^{p'}(Q(x,r)) \le c \gamma_a(Q(x,r)).$$
Thus taking $(p'-1)$th root and considering \eqref{i} if $x_i>r$ $i=1, \dots , n$, then
$$\gamma_a(Q(x,r))\le c \frac{x^ar^{|a|p+n-\nu p}}{\left(\prod_{i=1}^n\min\left\{1,\left(\frac{r}{x_i}\right)^{a_i}\right\}x^a\right)^{p}}\le c x^a r^{n-\nu p}.$$
If there is an $i$ such that $x_i\le r \le c_0$, then recalling that $Q(x,r)$ means the intersection of the cube with the positive orthant, $Q(x,r)\subset Q(y,r)$, where $x_i <y_i<x_i+r$ for all the suitable indices. Since the previous part is applicable to $Q(y,r)$, \eqref{v41} is proved.

To prove \eqref{v42} let $z\in Q\left(x,\frac{r}{2A\sqrt{n}}\right)$. If $\varrho>\frac{r}{A}$, $Q\left(x,\frac{r}{A\sqrt{n}}\right)\subset B(z,\varrho)$. Thus  as above we have
$$(G_{a,\nu}*_a\gamma|_{Q(x,r)})(z)\ge c \int_{\frac{r}{A}}^\infty  -g'(\varrho)\int_{Q(x,r)}T^u\chi_{Q\left(x,\frac{r}{A\sqrt{n}}\right)}(z)u^a d\gamma(u))d\varrho$$ $$\ge c g\left(\frac{r}{A}\right)\prod_{i=1}^n\min\left\{1,\left(\frac{r}{x_i}\right)^{a_i}\right\}\gamma_a(Q(x,r)).$$
Integrating this inequality over the cube $Q\left(x,\frac{r}{2A\sqrt{n}}\right)$ and the proceeding as above we get \eqref{v42}.

\medskip

\proof (of \eqref{Q} $\Rightarrow$ \eqref{V}) In view of Lemma \ref{l1}
$$G_{a,\nu}*_a(G_{a,\nu}*_a\gamma)^{p'}(x)=\int_0^\infty -g'(r)(\chi_{B(0,r)}*_a(G_{a,\nu}*_a\gamma)^{p'})(x)dr$$ $$=\int_0^1(\cdot)+\int_1^\infty(\cdot)=I+II,$$
where $g(r):=G_{a,\nu}(|x|)$, as above. With a fixed constant $m\ge 1$ we have
$$II\le c \left(\int_1^\infty -g'(r)\left(\chi_{B(0,r)}*_a(G_{a,\nu}*_a\gamma|_{Q(x,mr)})^{p'}\right)(x)dr\right.$$ $$\left.+\int_1^\infty -g'(r)\left(\chi_{B(0,r)}*_a(G_{a,\nu}*_a\gamma|_{Q(x,mr)^c})^{p'}\right)(x)dr\right)=c(II_1+II_2).$$
We decompose $I$ similarly, i.e.
$$I\le c(I_1+I_2).$$
According to Lemma \ref{tk} and by \eqref{Q} we have
$$II_1=\int_1^\infty -g'(r)\int_{B(x,r)}T^t\chi_{B(0,r)}(x)(G_{a,\nu}*_a\gamma|_{Q(x,mr)}(t))^{p'})t^adtdr$$ $$\le c \int_1^\infty -g'(r)\prod_{i=1}^n\min\left\{1,\left(\frac{r}{x_i}\right)^{a_i}\right\}\int_{Q(x,mr)}(G_{a,\nu}*_a\gamma|_{Q(x,mr)}(t))^{p'})t^adtdr$$ $$\le c \int_1^\infty -g'(r)\prod_{i=1}^n\min\left\{1,\left(\frac{r}{x_i}\right)^{a_i}\right\}\gamma_a\left(Q\left(x,\frac{r}{2\sqrt{n}}\right)\right)dr,$$
where the doubling property is applied. Thus applying Lemma \ref{tk} and the doubling property again we get
$$II_1\le c\int_1^\infty -g'(r)\int_{Q(x,\frac{r}{2\sqrt{n}})}T^t\chi_{Q(0,\frac{r}{\sqrt{n}})}(x)t^ad\gamma(t)dr$$ $$\le c\int_1^\infty -g'(r)\int_{B(x,r)}T^t\chi_{B(0,r)}(x)t^ad\gamma(t)dr\le c (G_{a,\nu}*_a\gamma)(x).$$
The estimation of $I_1$ is similar, only here the doubling property is not needed because for $r<1$ the kernel itself is doubling.
$$I_1=\int_0^1-g'(r)\left(\chi_{B(0,r)}*_a(G_{a,\nu}*_a\gamma|_{Q(x,mr)}(t))^{p'}\right)(x)dr$$ $$
\le c \int_0^1-g'(r)\prod_{i=1}^n\min\left\{1,\left(\frac{r}{x_i}\right)^{a_i}\right\}\int_{Q(x,mr)}(G_{a,\nu}*_a\gamma|_{Q(x,mr)}(t))^{p'})t^adtdr$$ $$ \le c \int_0^1-g'(r)\prod_{i=1}^n\min\left\{1,\left(\frac{r}{x_i}\right)^{a_i}\right\}\gamma_a(Q(x,mr))dr$$ $$\le c \int_0^1-g'(r)\int_{Q(0,2mr)}T^t\chi_{Q(0,4mr)}(x)\le c\int_0^1-g'(r)\chi_{Q(0,4mr)}*_a\gamma(x)dr$$ $$\le  c \int_0^1-g'(r)\chi_{B(0,4m\sqrt{n}r)}*_a\gamma(x)dr$$ $$=c\int_0^{4m\sqrt{n}}-g'(\varrho)\chi_{B(0,\varrho)}*_a\gamma(x)d\varrho\le c (G_{a,\nu}*_a\gamma)(x).$$
To estimate $I_2$ we need further decomposition.
$$I_2\le c \left(\int_0^1-g'(r)\left(\chi_{B(0,r)}*_a(G_{a,\nu}*_a\gamma|_{Q(x, 2m)\cap Q(x,mr)^c})^{p'}\right)(x)dr\right.$$ $$\left.+\int_0^1-g'(r)\left(\chi_{B(0,r)}*_a(G_{a,\nu}*_a\gamma|_{Q(x,2m)^c})^{p'}\right)(x)dr\right)=c(I_{2,1}+I_{2,2}).$$
Starting with $I_{2,2}$ we decompose the complementer set to cubes.
$$\mathbb{R}^{n}_+\setminus Q(x,2m)=\cup_{k=1}^\infty Q(y_k,1),$$
such that $y_{k,i}\ge 1$, that is the "first rows" of the cubes stand on the boundary hyperplanes and if $Q(x, 2m)\cap Q(y_k,1)\neq \emptyset$, then we take $\tilde{Q}(y_k,1)Q(y_k,1)\setminus Q(x, 2m)$.

Recalling that $t\in B(x,r)$, $r\le 1$, $u\in Q(y_k,1)$, $\|t-u\|>1$, in view of \eqref{v5}
$$D\le 6+\frac{2\|y_k-t\|}{\|y_k-x\|+\|u-t\|}\le 6+\frac{2\|y_k-x\|+1}{\|y_k-x\|+\|u-t\|}\le 6+2+\frac{1}{2}.$$ If $\vartheta \in [0,\pi)^n$ is such that $d(t,u,\vartheta)\ge d(x,y_k,\vartheta)$, substituting into $G_{a,\nu}$ the opposite inequality is valid. If $\vartheta \in [0,\pi)^n$ is such that $d(t,u,\vartheta)<d(x,y_k,\vartheta)$, than taking into consideration that $\|t-u\|>1$ i.e. $G_{a,\nu}$ behaves exponentially,
$$G_{a,\nu}(\sqrt{d(t,u,\vartheta)})=$$ \begin{equation}\label{vbe}G_{a,\nu}\left(\sqrt{d(x,y_k,\vartheta)}-\left(\sqrt{d(x,y_k,\vartheta)}-\sqrt{d(t,u,\vartheta)}\right)\right)\le c G_{a,\nu}(\sqrt{d(x,y_k,\vartheta)}).\end{equation}
Integrating the inequality above against $\sigma(\vartheta)$ we have
$$(G_{a,\nu}*_a\gamma|_{Q(x,2m)^c}(t))^{p'}=\left(\sum_{k=1}^\infty\int_{Q(y_k,1)}T^uG_{a,\nu}(t)u^ad\gamma(u)\right)^{p'}$$ $$\le c \left(\sum_{k=1}^\infty T^{y_k}G_{a,\nu}(x)\gamma_a(Q(y_k,1))\right)^{p'}$$ $$\le c\left(\sum_{k=1}^\infty T^{y_k}G_{a,\nu}(x)y_k^a\right)^{\frac{p'}{p}}\sum_{k=1}^\infty T^{y_k}G_{a,\nu}(x)y_k^{-a\frac{p'}{p}}\gamma_a(Q(y_k,1))^{p'}=:cAS,$$
where $A$ is a constant which is independent of $x$. Indeed,
$$\sum_{k=1}^\infty T^{y_k}G_{a,\nu}(x)y_k^a\le c\int_{\mathbb{R}^n_+}T^yG_{a,\nu}(x)y^ady = c\ G_{a,\nu}*_a\lambda(x)$$ $$\le c \ \int_0^\infty -g'(r)r^{|a|}\lambda(B(x,r))dr\le c \ \int_0^\infty -g'(r)r^{|a|+n}dr<\infty.$$
The last step follows from \eqref{i}, \eqref{j}, \eqref{k}.
Applying Lemma \ref{v4} to $S$ we have
$$S\le c \sum_{k=1}^\infty T^{y_k}G_{a,\nu}(x)y_k^{-a\frac{p'}{p}}\left(\prod_{i=1}^n\left(\max\{y_{k,i},1\}\right)^{a_i}\right)^{p'-1}\gamma_a(Q(y_k,1))$$ $$\le  c \sum_{k=1}^\infty T^{y_k}G_{a,\nu}(x)\gamma_a(Q(y_k,1)).$$
Thus, applying \eqref{v5} as above,
$$I_{2,2} \le c \int_0^1-g'(r)\int_{B(x,r)}T^t\chi_{B(0,r)}(x)\sum_{k=1}^\infty T^{y_k}G_{a,\nu}(x)\gamma_a(Q(y_k,1))t^adtdr$$ $$\le c\int_0^1-g'(r)\int_{B(x,r)}T^t\chi_{B(0,r)}(x)\sum_{k=1}^\infty \int_{Q(y_k,1)}T^{u}G_{a,\nu}(x)u^ad\gamma(u)t^adtdr$$ $$\le c\int_0^1-g'(r)r^{n+|a|}dr\left(G_{a,\nu}*_a\gamma|_{Q(x,2m)^c}\right)(x)\le c (G_{a,\nu}*_a\gamma)(x),$$
where \eqref{i} and \eqref{v9} were taken into consideration.

To estimate $I_{2,1}$ first we deal with convolution under the $p'$th power. Considering the doubling property by \eqref{v9} we have
$$\left(G_{a,\nu}*_a\gamma|_{Q(x,2m)\cap Q(x,mr)^c}\right)(t)\le c \left(G_{a,\nu}*_a\gamma|_{Q(x, 2m)\cap Q(x,mr)^c}\right)(x)$$ $$\le c \int_0^\infty -g'(\varrho)\left(\chi_{Q(0,\varrho)}*_a\gamma_{Q(x,2m)\setminus Q(x,mr)}\right)(x)d\varrho $$ $$=c\left(\int_{mr}^{2m}(\cdot)+\int_{2m}^\infty(\cdot)\right)=c(I_{2,1,1}+I_{2,1,2}).$$
By the doubling property and Lemma \ref{tk}
$$I_{2,1,1}\le c\int_{mr}^{2m}-g'(\varrho)\prod_{i=1}^n\min\left\{1,\left(\frac{\varrho}{x_i}\right)^{a_i}\right\}\int_{Q(x,\varrho)}u^ad\gamma(u)d\varrho$$ $$\le c\int_{mr}^{2m}-g'(\varrho)\gamma(Q(x,\varrho))x^a\prod_{i=1}^n\min\left\{1,\left(\frac{\varrho}{x_i}\right)^{a_i}\right\}d\varrho$$ $$\le c \int_{mr}^{2m}-g'(\varrho)\varrho^{|a|}\gamma(Q(x,\varrho))d\varrho.$$
Similarly
$$I_{2,1,2}\le \int_{2m}^\infty -g'(\varrho)\varrho^{|a|}\gamma(Q(x,2m))d\varrho\le c\gamma(Q(x,m))\int_{m}^{2m}-g'(\varrho)\varrho^{|a|}d\varrho.$$
Thus
\begin{equation}\label{i212}I_{2,1,1}+I_{2,1,2}\le c \int_{mr}^{2m}-g'(\varrho)\varrho^{|a|}\gamma(Q(x,\varrho))d\varrho,\end{equation}
and by \eqref{v9} and \eqref{i} we have
$$I_{2,1} \le c \int_0^1-g'(r)\int_{B(x,r)}T^t\chi_{B(0,r)}(x)t^adt\left(\int_{mr}^{2m}-g'(\varrho)\varrho^{|a|}\gamma(Q(x,\varrho))d\varrho\right)^{p'}dr$$ $$ \le c \int_0^1r^{\nu-1}\left(\int_{mr}^{2m}-g'(\varrho)\varrho^{|a|}\gamma(Q(x,\varrho))d\varrho\right)^{p'}dr$$ $$\le \int_0^{2m}\left(-g'(\varrho)\varrho^{|a|+1}\gamma(Q(x,\varrho))\right)^{p'}\varrho^{\nu-1}d\varrho,$$
where in the last step Hardy's inequality is applied. Thus, again by doubling property, cf. \eqref{v9} we have
$$I_{2,1}\le c \int_0^R\left(-g'(\varrho)\varrho^{|a|+1}\gamma(Q(x,\varrho))\right)^{p'}\varrho^{\nu-1}d\varrho+\int_R^{2m}(\cdot)$$
$$\le c \int_0^R\left(\int_\varrho^\infty -g'(t)t^{|a|}\gamma(Q(x,t)dt)\right)^{p'}\varrho^{\nu-1}d\varrho+ c\int_R^{2m}\varrho^{\nu(1-p)-1}d\varrho,$$
where Lemma \ref{v4} is applied. Thus
$$I_{2,1}\le c R^\nu(G_{a,\nu}*_a\gamma)^{p'}(x)+R^{-\nu(p-1)}.$$
If $B:=(G_{a,\nu}*_a\gamma)^{\frac{1-p'}{\nu}}(x)\le 2m$, then we choose $R=B$ and so
$$I_{2,1}\le c G_{a,\nu}*_a\gamma(x).$$
If $B>2m$, then
$$I_{2,1}\le c \int_0^{B}\left(-g'(\varrho)\varrho^{|a|+1}\gamma(Q(x,\varrho))\right)^{p'}\varrho^{\nu-1}d\varrho $$ $$\le c (G_{a,\nu}*_a\gamma)^{p'}(x)\int_0^{B}\varrho^{\nu-1}d\varrho \le c (G_{a,\nu}*_a\gamma)(x).$$

Before we turn to the estimation of $II_2$, we show that  $(G_{a,\nu}*_a\gamma)(x)$ is finite for almost all $x$.
Indeed,
$$\int_{Q(x,r)}\left((G_{a,\nu}*_a\gamma)(y)\right)^{p'}y^ady\le \int_{Q(x,r)}\left(\left(G_{a,\nu}*_a\gamma|_{Q(x,mr)}\right)(y)\right)^{p'}y^ady$$ $$+\int_{Q(x,r)}\left(\left(G_{a,\nu}*_a\gamma|_{Q(x,2m)\cap Q(x,mr)^c}\right)(y)\right)^{p'}y^ady$$ $$+\int_{Q(x,r)}\left(\left(G_{a,\nu}*_a\gamma|_{Q(x,2m)^c\cap Q(x,mr)^c}\right)(y)\right)^{p'}y^ady=J_1+J_2+J_3.$$
In view of \eqref{Q} we have
$$J_1\le \int_{Q(x,mr)}\left(\left(G_{a,\nu}*_a\gamma|_{Q(x,mr)}\right)(y)\right)^{p'}y^ady\le c\gamma_a(Q(x,mr)).$$
Assume for a moment that $r\ge 1$, say. If $r<1$, then the integral on $Q(x,r)$ is less than the integral on $Q(x,1)$. According to \eqref{i212} and Lemma \ref{v4} we get
$$J_2\le c\int_{Q(x,r)}\left(\int_{mr}^{2m} -g'(\varrho)\varrho^{|a|}\gamma(Q(x,\varrho))d\varrho\right)^{p'}y^ady$$ $$\le c\lambda_a(Q(x,r))x^{ap'} m^{|a|p'}\left(\int_{mr}^{2m}\varrho^{n-\nu p}\varrho^{|a|}\varrho^{-(n+|a|-\nu+1)}d\varrho\right)^{p'}$$ $$\le c(m,a,n,\nu, x,r)<\infty.$$
As we have seen in the estimation of $I_{2,2}$
$$J_3\le c\int_{Q(x,r)}\sum_{k=1}^\infty T^{y_k}G_{a,\nu}(x)\gamma_a(Q(y_k,1))y^ady$$ $$\le c\lambda_a(Q(x,r))\sum_{k=1}^\infty T^{y_k}G_{a,\nu}(x)y_k^a<\infty,$$
where Lemma \ref{v4} is applied, and we used the assumption that $y_{k,i}\ge 1$.
Thus the integral is finite which proves the statement.

Now we deal with $II_2$ recalling that $r\ge 1$, assuming that $m\ge 2$ and $(G_{a,\nu}*_a\gamma)(x)<\infty$.
$$\left(G_{a,\nu}*_a\gamma|_{Q(x,mr)^c}\right)(t)=\int_0^\infty -g'(\varrho)\int_{B(t,\varrho)\cap Q(x,mr)^c}T^u\chi_{B(0,\varrho)}(t)u^ad\gamma(u)d\varrho$$ $$\le c \int_{(m-1)r}^\infty -g'(\varrho)(\chi_{B(0,\varrho)}*_a\gamma)(t)d\varrho \le c\int_{(m-1)r}^\infty -g'(\varrho)\varrho^{|a|}\gamma(B(t,\varrho))d\varrho.$$
Since $\frac{\varrho +r}{\varrho}\le \frac{m}{m-1}\le 2$, by the doubling property we have
$$\left(G_{a,\nu}*_a\gamma|_{Q(x,mr)^c}\right)(t)\le c\int_{(m-1)r}^\infty -g'(\varrho)\varrho^{|a|}\gamma(B(t,\varrho+r))d\varrho $$ $$\le c\int_{(m-1)r}^\infty -g'(\varrho)\varrho^{|a|}\gamma(B(x,\varrho))d\varrho.$$
Considering that $g'\sim g$ on $\varrho>1$, in view of \eqref{i} $g'(\varrho)\le c e^{-\frac{\varrho}{p}}g'\left(\frac{\varrho}{p'}\right)$. Thus by H\" older's inequality
$$\left(\left(G_{a,\nu}*_a\gamma|_{Q(x,mr)^c}\right)(t)\right)^{p'}$$ $$\le c \int_{(m-1)r}^\infty -g'\left(\frac{\varrho}{p'}\right)^{p'}\varrho^{|a|}\gamma(B(x,\varrho))d\varrho \left(\int_{(m-1)r}^\infty e^{-\varrho}\varrho^{|a|}\gamma(B(x,\varrho))d\varrho\right)^{\frac{p'}{p}}.$$
By Lemma \ref{v4} and Remark \ref{rem2}
$$ \int_{(m-1)r}^\infty e^{-\varrho}\varrho^{|a|}\gamma(B(x,\varrho))d\varrho \le c\int_{(m-1)r}^\infty e^{-\varrho\left(1-\frac{p}{A}\right)}\varrho^{|a|(1-p)-n(p-1)}d\varrho,$$
so the second term is bounded by a constant which depends only on the parameters if $A=2p$, say. Thus, applying \eqref{i} again, as above but with $A=4p$, say, we have
$$II_2 \le c \int_1^\infty -g'(r)\int_{B(x,r)}T^t\chi_{B(0,r)}(x)\int_{(m-1)r}^\infty -g'(\varrho)\varrho^{|a|}\gamma(B(x,\varrho))d\varrho t ^adtdr$$
$$\le c \int_1^\infty -g'(r)r^{|a|}\gamma(B(x,r))\int_{(m-1)r}^\infty e^{-\varrho\left(\frac{1}{2}-\frac{p}{A}\right)}\varrho^{|a|(1-p)-n(p-1)}d\varrho$$ $$\le c \int_0^\infty  -g'(r) \chi_{B(0,r)}*_a\gamma(x)dr=c(G_{a,\nu}*_a\gamma)(x),$$
where the doubling condition is considered, and the proof is finished.

 \subsection{\eqref{V}$\Rightarrow$ \eqref{I}} To prove this implication we need the lemma below.

\medskip

\begin{lemma}\label{l210}
Let $f\ge 0$, $\mathrm{supp} f=:S$ is bounded ($S\subset B:=B(x_0,1)$, say for some $x_0$.) Assume that $f$ is uniformly bounded. Then there is a constant (independent of $\xi$) such that
$$\left(G_{a,\nu}*_a f(\xi)\right)^p\le c \ G_{a,\nu}*_a \left(f\left(G_{a,\nu}*_a f\right)^{p-1}\right)(\xi).$$

\end{lemma}

\medskip

For sake of selfcontainedness we cite the next lemma.

\begin{lemma}\label{vw}\cite[Lemma 4]{vw} Let $d$ be a quasi-metric, and $K(x,y)$ be a positive kernel function. Suppose $1<p<\infty$ and $\omega$ is a locally finite Borel measure on $X$. If $K(x,y)$ satisfies that
$$K(x,y)\le C_1 K(z,y), \ws \ws \mbox{if} \ws d(z,y)\le C_2 d(x,y),$$
and
$$K(x,y)\le C_1 K(x,z), \ws \ws \mbox{if} \ws d(z,x)\le C_2 d(x,y)$$
then
$$\left(\int_X K(x,y)d\omega(y)\right)^p\le C  \int_XK(x,y)\left(\int_XK(y,z)d\omega(z)\right)^{p-1}d\omega(y),$$
where $C$ depends on $C_1, p$ and the quasi-metric.

\end{lemma}

\medskip

\proof (of Lemma \ref{l210})
First we prove that
\begin{equation}\label{xiv} G_{a,\nu}*_a f(\xi)\le c \ G_{a,\nu}*_a f(v)\end{equation}
for all $v, \xi \in \mathbb{R}^{n}_+$, $v\in B$ and $\frac{3}{2}\le \|\xi-x_0\|$. By the assumption the expressions on both sides are finite.\\
First let $\|\xi-x_0\|>4\sqrt{n}+1$ and define $R:=\|\xi-x_0\|-1$. By Lemma \ref{tk} we have
$$G_{a,\nu}*_a f(\xi)=\int_R^\infty -g'(r)\int_{B(\xi,r)\cap S}T^\tau_{\chi(B(0,r)}(\xi)\tau^af(\tau)d\tau dr$$
\begin{equation}\label{el}\le c\int_R^\infty -g'(r)\int_{B(\xi,r)\cap S}\prod_{i=1}^n\min\left\{1, \left(\frac{r}{\xi_i}\right)^{a_i}\right\}\tau^af(\tau)d\tau dr.\end{equation}
On the other hand
$$G_{a,\nu}*_a f(v)\ge \int_{4\sqrt{n}}^\infty -g'(r)\int_{B(v,r)\cap S}T^t_{\chi(B(0,r)}(v)t^af(t)dt dr$$ $$\ge \int_{4\sqrt{n}}^\infty -g'(r)\int_{Q\left(v,\frac{r}{\sqrt{n}}\right)\cap S}Tt_{Q\left(v,\frac{r}{\sqrt{n}}\right)}(v)t^af(t)dt dr$$ \begin{equation}\label{ma}\ge c \int_{4\sqrt{n}}^\infty -g'(r)\int_{Q\left(v,\frac{r}{2\sqrt{n}}\right)\cap S}\prod_{i=1}^n\min\left\{1, \left(\frac{r}{v_i}\right)^{a_i}\right\}t^af(t)dt dr.\end{equation}
With the abbreviation $h_\xi(i):=\min\left\{1, \left(\frac{r}{\xi_i}\right\}\right)^{a_i}\tau_i^{a_i}$ we have if $\xi_i\ge v_i$, then $h_v(i)\ge h_\xi(i)$. If $\xi_i< v_i$ and $r<v_i$, then $h_v(i)\ge \left(r \left(1-\frac{r}{2\sqrt{n}v_i}\right)\right)^{a_i}>c \ r^{a_i}$ and $h_\xi(i)\le c \ r^{a_i}$, always. If $\xi_i<v_i<r$, then $h_\xi(i)=\tau_i^{a_i}$, $h_v(i)=t_i^{a_i}$, that is $h_\xi(i) \sim h_v(i)$, since $t, \tau \in S$. Considering that if $r>4\sqrt{n}$ then $Q\left(v,\frac{r}{2\sqrt{n}}\right)\cap S=S$ and comparing \eqref{el} and \eqref{ma} the proof of this case is finished.\\
If $\frac{3}{2}\le \|\xi-x_0\| \le 4\sqrt{n}+1$, we proceed as in \eqref{vbe}, i.e. in this case $T^\tau g(\xi)\le T^\tau g(v)$ which ensures \eqref{xiv}. Thus in these cases
$$\left(G_{a,\nu}*_a f(\xi)\right)^p=\int_BT^vG_{a,\nu}(\xi)\left(G_{a,\nu}*_a f(\xi)\right)^{p-1}v^af(v)dv$$ $$\le c \ \int_BT^vG_{a,\nu}(\xi)\left(G_{a,\nu}*_a f(v)\right)^{p-1}v^af(v)dv =c \ G_{a,\nu}*_a(f\left(G_{a,\nu}*_a f\right)^{p-1}(\xi).$$
Finally, let $\|\xi-x_0\| \le \frac{3}{2}$. Now we apply Lemma \ref{vw} on $B(x_0,1)$ with $d(x,y)=\|x-y\|$, $K(x,y)=T^y G_{a,\nu}(x)$ and $d\omega(y)=y^af(y)dy$. To this end it is enough to show, that if $x,y,z \in B(x_0,1)$ and $\|y-z\|\le B \|x-y\|$, $B>2$, then $T^y G_{a,\nu}(x)\le C_1 T^y G_{a,\nu}(z)$.\\
Let  $\|y-z\|\le B \|x-y\|$ and $\vartheta \in [0,\pi]^n$ such that $d(y,z,\vartheta)<1$. In view of \eqref{v51} there is a constant $C$ such that $d(y,z,\vartheta)\le C d(x,y,\vartheta)$. In this case $G_{a,\nu}$ behaves as a power function (with negative exponent), that is  $G_{a,\nu}(\sqrt{d(x,y,\vartheta})\le c_1 G_{a,\nu}(\sqrt{d(y,z,\vartheta})$. If $\vartheta \in [0,\pi]^n$ such that $d(y,z,\vartheta)\ge 1$, we apply \eqref{v5}, and as above (see \eqref{vbe} again) we have that $G_{a,\nu}(\sqrt{d(x,y,\vartheta})\le c_2 G_{a,\nu}(\sqrt{d(y,z,\vartheta})$. Thus integrating these inequalities on $[0,\pi]^n$ against the probability measure $\sigma$, we get the required inequality and the proof of this part is finished by applying Lemma \ref{vw}.

\medskip

\proof(of \eqref{V} $\Rightarrow$ \eqref{I}) First observe that it is enough to prove the statement for nonnegative uniformly bounded functions, and we can assume that the diameter of the support of the function is at most one. Indeed, as usually, $f$ can be decomposed to positive and negative parts. Then $\mathbb{R}^n_+$ is decomposed to cubes ($Q_i$) of diameter at most one. Assuming that $f\ge 0$, let $f_{N,i}=\min\{f|_{Q_i},N\}$ and $f_N= \min \{f,N\}$. If $f_{N,i}$ fulfils \eqref{I} for all $i$ then the linearity of convolution and the triangle inequality ensures that $f_N$ fulfils \eqref{I}. Finally applying Fatou's lemma the inequality can be derived to $f$. Thus below we can apply Lemma \ref{l210}.
\begin{equation}\label{bb}\|G_{a,\nu}*_a f\|_{p,\gamma_a}^p\le c \ \int_{\mathbb{R}^n_+}G_{a,\nu}*_a \left(f\left(G_{a,\nu}*_a f\right)^{p-1}\right)(\xi)\xi^ad\gamma(\xi)\end{equation} $$=\int_{\mathbb{R}^n_+}\int_{\mathbb{R}^n_+}T^tG_{a,\nu}(\xi)\xi^ad\gamma(\xi)f(t)\left(G_{a,\nu}*_a f(t)\right)^{p-1}t^adt\le c \ \|f\|_{p,a}\|(G_{a,\nu}*_a f\|_{p,\nu_a}^{p-1},$$
where $d\nu(t)=\left(G_{a,\nu}*_a\gamma(t)\right)^{p'}dt$. Repeating the same track of thoughts with $\nu$ instead of $\gamma$ we arrive to
\begin{equation}\label{ff}\|G_{a,\nu}*_a f\|_{p,\nu_a}^p\le c \ \|f\|_{p,a}\|(G_{a,\nu}*_a f\|_{p,\nu_{1,a}}^{p-1},\end{equation}
where
$$d\nu_1(t)=\left(G_{a,\nu}*_a\nu(t)\right)^{p'}dt=\left(G_{a,\nu}*_a\left(G_{a,\nu}*_a\gamma(t)\right)^{p'}\right)^{p'}dt\le c \ \left(G_{a,\nu}*_a\gamma(t)\right)^{p'}dt.$$
Thus if $\|G_{a,\nu}*_a f\|_{p,\nu_a}=:M$ is finite, equation \eqref{ff} can be divided by $M^{p-1}$ and comparing it with equation \eqref{bb} we have \eqref{I}.

To show that $M<\infty$, first we observe that $G_{a,\nu}\in L^1_a$ so in view of \eqref{e} $G_{a,\nu}*_a f(x)<c \ \|f\|_\infty$. Let $S:=\mathrm{supp}f$.
$$\|G_{a,\nu}*_a f\|_{p,\nu_a}^p=\int_{ \mathbb{R}^{n}_+}\left(G_{a,\nu}*_a f(x)\right)^p\left(G_{a,\nu}*_a \gamma(x)\right)^{p'}x^adx$$ $$ \le\|G_{a,\nu}*_a f\|_\infty^{p-1}\int_{ \mathbb{R}^{n}_+}\left(G_{a,\nu}*_a f(x)\right)\left(G_{a,\nu}*_a \gamma(x)\right)^{p'}x^adx$$ $$=\|G_{a,\nu}*_a f\|_\infty^{p-1}\int_{ \mathbb{R}^{n}_+}\left(G_{a,\nu}*_a\left(G_{a,\nu}*_a \gamma\right)^{p'}\right)(t)f(t)t^adt$$ $$\le c \ \|f\|_\infty^p \int_{S}\left(G_{a,\nu}*_a\left(G_{a,\nu}*_a \gamma\right)^{p'}\right)(t)t^adt$$ $$\le c \ \|f\|_\infty^p \int_{S}G_{a,\nu}*_a \gamma(t)t^adt.$$
Now we show that
\begin{equation}\label{fi}\int_{S}G_{a,\nu}*_a \gamma(t)t^adt<\infty.\end{equation}
To this we prove that $G_{a,\nu}*_a\gamma \in L_{\mathrm{loc},a}^{p'}$. Then, considering that $\mathrm{diam} \ S\le 1$, H\" older inequality ensures \eqref{fi}.

Let us introduce the next notation (see \cite{h})
\begin{equation}\label{Md}M_{a,\nu,1}\gamma(x):=\sup_{0<r\le 1}\frac{\chi_{B_+(0,r)}*_a\gamma(x)}{r^{n+|a|-\nu}}.\end{equation}
Observe that
$$G_{a,\nu}*_a\gamma(x)\ge \int_r^\infty-g'(\varrho)\chi_{B(0,\varrho)}*_a\gamma(x)d\varrho \ge \chi_{B(0,r)}*_a\gamma(x)g(r).$$
Thus for all $r\le 1$ $G_{a,\nu}*_a\gamma(x)\ge \frac{c}{r^{n+|a|-\nu}} \chi_{B(0,r)}*_a\gamma(x)$, which implies that
\begin{equation}G_{a,\nu}*_a\gamma(x)\ge c M_{a,\nu,1}\gamma(x).\end{equation}
Applying this in view of \eqref{V} for almost all $x$ we have
$$ M_{a,\nu,1}\left(G_{a,\nu}*_a\gamma\right)^{p'}(x)\le c G_{a,\nu}*_a\left(G_{a,\nu}*_a\gamma\right)^{p'}(x)\le c \ (G_{a,\nu}*_a\gamma)(x)<\infty.$$
On the other hand for all $r\le 1$
$$ M_{a,\nu,1}\left(G_{a,\nu}*_a\gamma\right)^{p'}(x)=\sup_{0<r\le 1}\frac{1}{r^{n+|a|-\nu}}\int_{B(x,r)}T^t\chi_{B(0,r)}(x)\left(G_{a,\nu}*_a\gamma\right)^{p'}(t)t^adt$$ $$\ge \frac{1}{r^{n+|a|-\nu}}\int_{Q\left(x,\frac{r}{\sqrt{n}}\right)}T^t\chi_{Q\left(0,\frac{r}{\sqrt{n}}\right)}(x)\left(G_{a,\nu}*_a\gamma\right)^{p'}(t)t^adt$$ $$ \ge c \frac{1}{r^{n+|a|-\nu}}\int_{Q\left(x,\frac{r}{2\sqrt{n}}\right)}\prod_{i=1}^n\min\left\{1,\frac{r}{x_i\sqrt{n}}\right\}^{a_i}\left(G_{a,\nu}*_a\gamma\right)^{p'}(t)t^adt.$$
If we choose $r$ such that $r<\min\{1,x_i\sqrt{n}\}$ for all $i=1,\dots ,n$,\\ then $\frac{rt_i}{x_i\sqrt{n}}\ge \frac{r}{x_i\sqrt{n}}(x_i-\frac{r}{2\sqrt{n}})\ge \frac{r}{2\sqrt{n}}$. Thus
$$ c(n,a)r^{n-\nu} x^a(G_{a,\nu}*_a\gamma)(x)\ge\int_{Q\left(x,\frac{r}{2\sqrt{n}}\right)}\left(G_{a,\nu}*_a\gamma\right)^{p'}(t)t^adt.$$

\section{Schr\"odinger-Bessel operators}

Schr\"odinger operator on the halfline, say, with the potential $q_m:=\left(m^2-\frac{1}{4}\right)\frac{1}{x^2}$ is a Bessel operator, see \cite{le}. As it has widespread applications, it is examined by several authors. We mention only a few examples. For instance its domain is investigated under different extensions. Assuming that $m$ is complex, see \cite{dg}, and translated the singularity to $a\in\mathbb{R}$, see \cite{gps}. The number of negative eigenvalues is estimated in \cite{bs}. \\
Here the expression "Schr\"odinger-Bessel operator" stands for  
$$-\Delta_a-v=-\sum_{i=1}^n\left(\frac{\partial^2}{\partial x_i^2}+\frac{2\alpha_i+1}{x_i}\frac{\partial}{\partial x_i}\right) -v,$$ 
cf. \eqref{sre}. We deal with the least eigenvalue problem. The "method of cubes" appeared in \cite{f} with respect to classical Schr\"odinger operator. This result is improved by trace inequality in \cite{ks}. We extend this method to the Schr\"odinger-Bessel case.

\medskip

In this section $L^2_a(\mathbb{R}^n_+)$ means the space of all measurable functions which are even with respect to each variable $x_i$ $i=1, \dots , n$ and $\|(\cdot)\|_{2,a}$ is bounded.

\medskip

\note
$$\mathcal{S}_{ev}:=\left\{f\in C_{ev}^\infty(\overline{\mathbb{R}^n}_+) :\sup_{x\in \mathbb{R}^n_+}\left|x^\delta D^\beta f\right|<\infty \ws \ws \forall \ \delta,\beta \in \mathbb{N}^n\right\},$$
where $C_{ev}^\infty(\overline{\mathbb{R}^n}_+)=\cap C_{ev}^m(\overline{\mathbb{R}^n}_+)$, where the intersection taken for all finite $m$ and
$$C_{ev}^m(\overline{\mathbb{R}^n}_+):=\left\{f\in C^m(\overline{\mathbb{R}^n}_+ : \frac{\partial^{2k+1}f}{\partial x_i^{2k+1}}|_{x_i=0}=0 \ws \forall \ 0\le k \le \frac{m-1}{2}, \ws i=1, \dots, n\right\}.$$

\begin{remark}
 Let us denote by $G_{a,\nu}^c(x)=G_{a,\nu}(cx)$ and $V^c(x)=c^{n+|a|}v(cx)$ and let in Theorem \ref{vm} the measure $d\gamma(x)=v(x)dx$. By a simple replacement we have that \eqref{I}, \eqref{Q} and \eqref{V} fulfils with $G_{a,\nu}$, $v$ and $A_i$, respectively if and only if it fulfils with $G_{a,\nu}^c$, $V^c$ and the same  $A_i$. Thus Theorem \ref{vm} fulfils with the generalized Bessel kernel $G_{a,\nu}^c$ as well.
 \end{remark}

\medskip

\proof (of Theorem \ref{eiv})
Let $\mathcal{D}(H_a)\subset L^2_a(\mathbb{R}^n_+)$ is the domain of definition of $H_a$. $\mathcal{S}_{ev}$ is dense in $\mathcal{D}(H_a)$. Thus, by Rayleigh-Ritz method we have
$$-\lambda_1=-\inf_{u\in \mathcal{D}(H_a)}\frac{\langle H_a u,u\rangle_a}{\langle u,u\rangle_a}=\sup_{u\in \mathcal{S}_{ev}}\frac{\int_{\mathbb{R}^n_+}\left(|u(x)|^2v(x)-u(x)(-\Delta_a) u(x)\right)x^adx}{\int_{\mathbb{R}^n_+}|u(x)|^2x^adx}$$ $$=\inf\left\{\beta>0:\int_{\mathbb{R}^n_+}|u(x)|^2v(x)x^adx\le \int_{\mathbb{R}^n_+}\left(u(x)(-\Delta_a) u(x)+\beta|u(x)|^2\right)x^adx\right\}.$$
According to \eqref{hb} and \eqref{hn}
$$\langle (-\Delta_a)u+\beta u, u\rangle_a=k(n,a)\left\|\sqrt{(|\cdot|^2+\beta)}\hat{u}\right\|_{2,a}^2,$$
where $k(n,a)=\frac{2^{n-|a|}}{\prod_{i=1}^n\Gamma^2(\alpha_i+1)}$. Take now $\sqrt{(|\cdot|^2+\beta)}\hat{u}=:\hat{f}$. Considering \eqref{hG} $$\mathcal{H}_a^{-1}\frac{1}{\sqrt{(|\xi|^2+\beta)}}=\beta^{\frac{n-1+|a|}{2}}G_{a,1}\left(\sqrt{\beta}x\right)=:G_{a,1}^\beta(x),$$ thus by \eqref{hkonv} and \eqref{hn} again we have
$$-\lambda_1=\inf\left\{\beta>0:\int_{\mathbb{R}^n_+}(G_{a,1}^\beta*_af)^2(x)x^av(x)dx\le \|f\|_{2,a}^2, \ws f\ge 0\right\}.$$
$C_{\beta}$ stands for the least constant in
$$\int_{\mathbb{R}^n_+}(G_{a,1}^\beta*_af)^2(x)x^av(x)dx\le C_\beta\|f\|_{2,a}^2 \ws \ws \forall \ f\ge 0.$$
($-\lambda_1=\inf\{\beta : C_\beta\le 1\}$.) In view of Theorem \ref{vm}
$$C_\beta \sim \sup_Q\frac{1}{v_a(Q)}\int_Q (G_{a,1}^\beta*_av|_Q(x))^{2}x^adx=:M_0.$$
We show that the supremum above need only be taken over cubes with edge lengths ($l(Q)$) at most $2\beta^{-\frac{1}{2}}$. To this end, set $$M:=\sup_{Q, \ l(Q)\le \beta^{-\frac{1}{2}}}\frac{1}{v_a(Q)}\int_{\mathbb{R}^n_+} (G_{a,1}^\beta*_av|_Q(x))^{2}x^adx.$$
Let $Q$ be a cube such that $l(Q)>\beta^{-\frac{1}{2}}$. Then we decompose $Q=\cup_{i=1}^LQ_i$ such that $Q_i=Q\left(y_i,l\right)$, where $\frac{1}{2}\beta^{-\frac{1}{2}}\le l\le \beta^{-\frac{1}{2}}$, $\mathrm{int} Q_i\cap \mathrm{int} Q_j=\emptyset$ for all $1\le i,j\le L$ (L depends on $Q$). We have
$$\int_{\mathbb{R}^n_+} (G_{a,1}^\beta*_av|_Q(x))^{2}x^adx=\sum_{i,j}\int_{\mathbb{R}^n_+} (G_{a,1}^\beta*_av|_{Q_i})(x)(G_{a,1}^\beta*_av|_{Q_j})(x)x^adx$$ $$=\sum_{i=1}^L\sum_{k=0}^{L-1}\sum_{j, \mathrm{dist}(Q_i,Q_j)=kl}\int_{\mathbb{R}^n_+} (G_{a,1}^\beta*_av|_{Q_i})(x)(G_{a,1}^\beta*_av|_{Q_j})(x)x^adx$$
If $k=0$ or $k=1$, we take
\begin{equation}\label{sro1}(G_{a,1}^\beta*_av|_{Q_i})(G_{a,1}^\beta*_av|_{Q_j})\le \frac{1}{2}\left((G_{a,1}^\beta*_av|_{Q_i})^2+(G_{a,1}^\beta*_av|_{Q_j})^2\right).\end{equation}
If $k>1$, then by \eqref{hn} and \eqref{j} we have
$$I_{ij}:=\int_{\mathbb{R}^n_+} (G_{a,1}^\beta*_av|_{Q_i})(x)(G_{a,1}^\beta*_av|_{Q_j})(x)x^adx=\int_{Q_j}(G_{a,2}^\beta*_av|_{Q_i})(x)v(x)x^adx$$ $$=\int_{Q_j}\int_{Q_i}T^tG_{a,2}^\beta(x)t^av(t)x^av(x)dtdx$$ $$\le c\int_{Q_j}\int_{Q_i}\int_{[0,\pi]^n}e^{-\sqrt{\beta\left(\|x-t\|^2+4\sum_{k=1}^nx_kt_k\sin^2\frac{\vartheta_k}{2}\right)}}(tx)^\frac{a}{2}d\sigma(\vartheta)(tx)^\frac{a}{2}v(t)v(x)dtdx.$$
Recalling that the distance of $Q_i$ and $Q_j$ is at least $kl$, by convexity we have
$$J:=\int_{[0,\pi]^n}(\cdot)\le c(a) e^{-kc}\prod_{k=1}^n\int_0^\pi e^{-\sqrt{\frac{2\beta}{n}}\sqrt{x_kt_k}\sin\frac{\vartheta_k}{2}}\sqrt{x_kt_k}^{a_k}\sin^{a_k-1}\vartheta_k d\vartheta_k.$$
To show that the integral in the product is bounded by a constant which depends only on $a$ and $n$ we investigate the integral against $\vartheta_k$ and for sake of simplicity we introduce some abbreviations.
$$\int_0^\pi e^{-cb_k\sin\frac{\vartheta_k}{2}}b_k^{a_k}\sin^{a_k-1}\vartheta_k d\vartheta_k=\int_0^{\frac{\pi}{2}}(\cdot)+\int_{\frac{\pi}{2}}^\pi(\cdot)=I_k+II_k.$$
Since in the second integral $\sin\frac{\vartheta_k}{2}\ge \frac{1}{\sqrt{2}}$, this term is bounded by a constant which depends only on $a_k$ and $n$.
Taking into consideration that on $\left[0,\frac{\pi}{2}\right]$ $\sin x\sim x$, the first integral is estimated as $I_A=\int_0^{A} e^{-cbx}b^a x^{a-1}dx$. After integration by parts $[a]$ times we have
$$I_A \le e^{-bAc}\left(\sum_{m=1}^{[a]}-d(m,a)b^{a-m}\frac{A^{a-m}}{c^m}\right)+d([a],a)\frac{b^{a-[a]}}{c^{[a]}}\int_0^Ae^{-cbx}x^{a-[a]-1}dx.$$
Replacing in the last integral $u=bx$ we have\\ $\int_0^Ae^{-cbx}x^{a-[a]-1}dx\le b^{[a]-a}\int_0^\infty e^{-cu}u^{a-[a]-1}du\le c(a,n) b^{[a]-a}$, and so
$I_k\sim I_A \le c(a_k,n)$, and
$$J\le c(a,n) e^{-kc}.$$
Thus we have
\begin{equation}\label{sro2}I_{ij}\le c(n,a) e^{-kc}v_{\frac{a}{2}}(Q_i)v_{\frac{a}{2}}(Q_j).\end{equation}
If $u, t\in Q_i$, then
$$T^uG_{a,1}^\beta(t)$$ $$\ge c(\beta)\int_{\vartheta \in [0,\pi]^n, \atop \sin\frac{\vartheta_k}{2}\le \frac{1}{2\sqrt{n}\max\{l,y_{i,k}\}}}g\left(\sqrt{\beta \left(\|u-t\|^2+4\sum_{k=1}^nu_kt_k\sin^2\frac{\vartheta_k}{2}\right)}\right)d\sigma(\vartheta)$$ $$\ge c(a,n,\beta)\frac{1}{\prod_{k=1}^n\left(\max\{l,y_{i,k}\}\right)^{a_k}}g(\sqrt{\beta(\|u-t\|^2+1)})$$ $$=c(a,n,\beta)\frac{1}{\prod_{k=1}^n\left(\max\{l,y_{i,k}\}\right)^{a_k}}.$$
In view of Remark \ref{rem2}
$$(G_{a,1}^\beta*_a\chi|_{Q_i})(t)=\int_{Q_i}T^uG_{a,1}^\beta(t)u^adu\ge c \lambda(Q_i).$$

So by Fubini's theorem
$$v_a(Q_i)\le \frac{c}{\lambda(Q_i)}\int_{Q_i}(G_{a,1}^\beta*_a\chi|_{Q_i})(t)t^av(t)dt=\frac{c}{\lambda(Q_i)}\int_{Q_i}(G_{a,1}^\beta*_av|_{Q_i})(u)u^adu.$$
Since $v$ is doubling, as above we have
$$2v_{\frac{a}{2}}(Q_i)v_{\frac{a}{2}}(Q_j)\le \frac{v_a^2(Q_i)}{\prod_{k=1}^n\left(\max\{l,y_{i,k}\}\right)^{a_k}}+\frac{v_a^2(Q_j)}{\prod_{k=1}^n\left(\max\{l,y_{j,k}\}\right)^{a_k}}$$ $$\le  \frac{c}{\lambda(Q_i)^2}\left(\frac{1}{\prod_{k=1}^n\left(\max\{l,y_{i,k}\}\right)^{a_k}}\left(\int_{Q_i}(G_{a,1}^\beta*_av|_{Q_i})(u)u^adu\right)^2\right.$$ $$\left.+\frac{1}{\prod_{k=1}^n\left(\max\{l,y_{j,k}\}\right)^{a_k}}\left(\int_{Q_j}(G_{a,1}^\beta*_av|_{Q_j})(u)u^adu\right)^2\right)$$
\begin{equation}\label{sro3}\le c(l) \left(\int_{Q_i}\left((G_{a,1}^\beta*_av|_{Q_i})(u)\right)^2u^adu+\int_{Q_j}\left((G_{a,1}^\beta*_av|_{Q_j})(u)\right)^2u^adu\right).\end{equation}
Notice that for a fixed cube $Q_i$ there are at most $ck^{n-1}$ cubes $Q_j$ at a distance $kl$, in view of \eqref{sro1}, \eqref{sro2} and \eqref{sro3}
$$\sum_{i,j}\int_Q (G_{a,1}^\beta*_av|_{Q_i})(x)(G_{a,1}^\beta*_av|_{Q_j})(x)x^adx $$ $$\le c\left(2+\sum_{k=2}^\infty e^{-k}k^{n-1} \right)\sum_{i=1}^L \int_{Q_i}(G_{a,1}^\beta*_av|_{Q_i})^2(u)u^adu.$$
Finally we have
$$\int_{\mathbb{R}^n_+} (G_{a,1}^\beta*_av|_Q(x))^{2}x^adx$$ $$\le c\sum_{i=1}^L\int_{Q_i}(G_{a,1}^\beta*_av|_{Q_i})^2(u)u^adu\le cM\sum_{i=1}^Lv_a(Q_i)=cMv_a(Q)$$
for each cube $Q$. Taking supremum we have
$$M_0\sim M.$$
Applying the Fourier-Bessel transformation again we have $\int_{\mathbb{R}^n_+} (G_{a,1}^\beta*_av|_Q(x))^{2}x^adx=\int_{\mathbb{R}^n_+} (G_{a,2}^\beta*_av|_Q(x))v_Q(x)x^adx$. We finish the proof with the observation (again on the transformed side) that $$\int_{\mathbb{R}^n_+} (G_{a,2}^\beta*_av|_Q(x))v_Q(x)x^adx\sim \int_{\mathbb{R}^n_+} (G_{a,2}*_av|_Q(x))v_Q(x)x^adx,$$
which after a renormalization leads to the inequalities of Theorem \ref{eiv}.

\vspace{5mm}

{\small{\noindent Both authors:\\
Department of Analysis and Operations Research,\\ Institute of Mathematics,\newline
Budapest University of Technology and Economics \newline
 M\H uegyetem rkp. 3., H-1111 Budapest, Hungary.

 \vspace{3mm}

\noindent \'A. P. Horv\'ath: g.horvath.agota@renyi.hu\\
M. Chegaar: mouna.chegaar@edu.bme.hu}

\end{document}